\documentclass[12pt]{article}
\usepackage{amsmath,amsthm,amsfonts,amscd,eucal}
\begin{document}
\def\RR{{\mathbb R}}
\def\CC{{\mathbb C}}
\def\ZZ{{\mathbb Z}}
\def\s1{{S^1}}
\def\sl2{{{\rm SL}(2,\RR)}}
\def\psl2{{{\rm PSL}(2,\RR)}}
\def\u1{{\rm U(1)}}
\def\A{{\cal A}}
\def\B{{\cal B}}
\def\C{{\cal C}}
\def\D{{\cal D}}
\def\F{{\cal F}}
\def\H{{\cal H}}
\def\I{{\cal I}}
\def\K{{\cal K}}
\def\k{{\rm K}}
\def\M{{\cal M}}
\def\N{{\cal N}}
\def\O{{\cal O}}
\def\P{{\cal P}}
\def\R{{\cal R}}
\def\T{{\cal T}}
\def\U{{\cal U}}
\def\V{{\cal V}}
\def\W{{\cal W}}
\def\G{{\bf G}}

\newtheorem{theorem}{Theorem}[section]
\newtheorem{definition}[theorem]{Definition}
\newtheorem{corollary}[theorem]{Corollary}
\newtheorem{proposition}[theorem]{Proposition}
\newtheorem{lemma}[theorem]{Lemma}
\newtheorem{remark}[theorem]{Remark}
\title{\bf The Virasoro Algebra and Sectors with
Infinite Statistical Dimension}

\author{SEBASTIANO CARPI\footnote{Supported in part by the Italian MIUR 
and GNAMPA-INDAM}
\\ 
 Dipartimento di Scienze \\
Universit\`a ``G. d'Annunzio'' di Chieti-Pescara\\
Viale Pindaro 42, 65127 Pescara -- Italy\\
e.mail: carpi@sci.unich.it }
\date{}
\maketitle
\renewcommand{\sectionmark}[1]{}
\begin{abstract} 
We show  that the sectors with lowest weight $h\geq
0$, $h\neq j^2$, $j\in \frac{1}{2}\ZZ$ of the local net of von
Neumann algebras on the circle generated by the Virasoro algebra with
central charge $c=1$ have infinite statistical dimension.
\end{abstract}

\newpage

\section{Introduction}
The notion of statistical dimension of superselection sectors,
introduced by Doplicher, Haag, and Roberts in \cite{DHR} is one of the most 
important concepts emerging in the formulation of 
Quantum Field Theory in therms of local nets of
operator algebras ( see \cite{Haag} for a general reference on this
subject). The  deep connection with Jones' theory on index
for subfactors \cite{jones,kosaki}, established by Longo 
\cite{longo 89-90} is a remarkable illustration of the relevance of this
notion. 

For an irreducible  representation $\pi$ of the algebra of observables
$\A$ satisfying the DHR selection criterion
the finiteness of the (statistical) dimension $d(\pi)$ is equivalent to
the existence of a conjugate representation $\overline{\pi}$ corresponding
to the particle-antiparticle symmetry \cite{DHR}, a condition which is
very natural on physical grounds. In fact for local nets over a 
four dimensional Minkowski space-time no example of (irreducible) sector 
with infinite dimension is known and the possibility  that in this context
the existence of such sectors can be excluded for physically reasonable
algebras of observables is  still open.

The situation is different in the case of
conformal nets on $\s1$, i.e. nets associated to chiral components of
2D conformal field theories, where irreducible
representations with infinite dimension seem to appear naturally.
Examples have been found by Fredenhagen \cite{fredenhagen} and 
Rehren has given arguments indicating that for the nets generated by the
Virasoro with central charge $c\geq 1$ most of the irreducible representations 
should have infinite dimension \cite{Reh}. Moreover the analysis of these 
representations in a model independent framework has been initiated in
\cite{BCL}  

In this note we show (Theorem \ref{infcont}), in agreement with
the arguments in \cite{Reh}, that the representations
of the Virasoro algebra with central charge $c=1$ and lowest weight $h\geq 
0$, $h\neq j^2$, $j\in \frac{1}{2}\ZZ$
give rise to representations with infinite dimension of the corresponding
conformal net $\A_{\rm Vir}$.

Our strategy of proof differs from the one adopted in \cite{fredenhagen}
where (partial) computation of fusion rules is used to infer infinite
dimension. Part of the fusion rules for the Virasoro algebra with
$c=1$ have been recently computed by Rehren and Tuneke \cite{RT} but we 
shall not use their results. 

Instead of the fusion structure we use a formula, which appeared in 
\cite{Reh2}, giving the dimension the
restriction of a representation of a net $\A$ to a subsystem $\B\subset \A$
(Proposition \ref{indrest} in this note) and well known results on the
representation theory of the Virasoro algebra \cite{Kac}. 
As another interesting
application of this formula we show, generalizing a result in 
\cite{xu 2000b}, that for finite index subsystems of certain rational nets
twisted sectors always exist (Proposition \ref{twisted}).

  \section{Conformal nets, their representations and subsystems}
Let $\I$ be the set of nonempty, nondense,  open  intervals of
unit circle $\s1$.

A {\it conformal net on $\s1$} is a family
$\A=\{ \A(I) | I \in \I \}$
of von Neumann algebras, acting on a infinite-dimensional separable
Hilbert space ${\H}_\A$, satisfying the following properties: 
\begin{itemize}
\item[(i)] {\it Isotony.}
\begin{equation}
\A(I_1) \subset \A (I_2)\;\;{\rm for} \;I_1 \subset I_2,
\;\;
I_1, I_2 \in \I.
\end{equation}

\item[(ii)] {\it Locality.}
\begin{equation}
\A(I_1) \subset \A(I_2)^{\prime}\;\;{\rm for }\;
I_1 \cap I_2 =\emptyset ,\;I_1, I_2 \in \I.
\end{equation}

\item[(iii)] {\it Conformal covariance.} There exists  a strongly
continuous unitary representation $U$ of $\psl2$ in ${\H}_\A$ such that
\begin{equation}
U({\alpha}){\cal A}(I)U({\alpha})^{-1}={\cal A}({\alpha}I) \;\;{\rm for}\;I
\in \I,\; \alpha\in \psl2,
\end{equation}

where $\psl2$ acts on $\s1$ by Moebius transformations. 
 
\item[(iv)] {\it Positivity of the energy.} The {\it conformal
Hamiltonian}
$L_{0}$,
which generates the restriction of $U$
to the one-parameter group of rotations
has non negative spectrum.

\item[(v)] {\it Existence of the vacuum.} There exists a unique (up to a
phase)  $U$-invariant unit vector $\Omega\in \H_\A$.

\item[(vi)] {\it Cyclicity of the vacuum.} ${\Omega}$ is cyclic for the algebra
$\A(\s1):=\bigvee_{I\in \I} \A (I)$ 
\end{itemize} 

Some consequences of the axioms are \cite{FrG,GuLo96}:
\begin{itemize}
\item[(vii)] {\it Reeh-Schlieder property.} For every $I\in \I$,
$\Omega$ is cyclic and separating for $\A(I)$.

\item[(viii)] {\it Haag duality.} For every $I \in \I$
\begin{equation}
\A(I)^{\prime}= \A(I^{c}),
\end{equation}
where $I^{c}$ denotes the interior of $S^{1}\backslash I$.

\item[(ix)] {\it Factoriality.} The algebras ${\A}(I)$ are type
${\rm III}_{1}$ factors.
\end{itemize}

A conformal net $\A$ is said to be {\it split} if given two intervals
$I_1,I_2\in\I$ with the closure of $I_1$ contained in $I_2$, there exists
a type I factor $\N(I_1,I_2)$ such that 
\begin{equation}
\A(I_1)\subset \N(I_1,I_2) \subset \A(I_2).
\end{equation}
Moreover, if for every $I,I_1,I_2\in\I$ with $I_1,I_2$ obtained
by removing a point from $I$ we have
\begin{equation}
\A(I_1)\vee\A(I_2)=\A(I),
\end{equation} 
then $\A$ is said to be {\it strongly additive}. The split property and
strong additivity do not follow from the axioms of conformal nets but they  
are satisfied in many interesting models.
 
A {\it representation} of a conformal net $\A$ is a family 
$\pi =\{ \pi_I| \; I\in \I \}$ where $\pi_I$ is a representation of $\A(I)$ on
a fixed Hilbert space $\H_\pi$, such that 

\begin{equation}
\pi_J |_{\A(I)} =\pi_I \;\; {\rm for} \;I\subset J.
\end{equation}
 
Irreducibility, direct sums and unitary equivalence of representations
of conformal nets can be defined in a natural way, see \cite{FrG,GuLo96}. 
The unitary equivalence class of an irreducible representation $\pi$ on 
a separable Hilbert space is called a {\it sector} and denoted $[\pi]$. 
The identical representation of $\A$ on $\H_\A$ is called the {\it vacuum
representation} and it is irreducible. The corresponding sector is called 
the {\it vacuum sector}. 

If $\H_\pi$ is separable then $\pi$ is automatically locally normal, namely 
$\pi_I$ is normal for each $I\in \I$ and hence $\pi_I(\A(I))$ is a type 
${\rm III}_1$ factor.
A representation $\pi$ is said to be {\it covariant} if there is a
strongly
continuous unitary representation $U_\pi$ on $\H_\pi$ of the universal
covering group $\widetilde{\psl2}$ of $\psl2$ such that 

\begin{equation}
\label{covariance}
{\rm Ad} U_\pi(\alpha)\pi_I=\pi_{\alpha I}{\rm Ad} U(\alpha),
\end{equation}
where $U$ has been lifted to $\widetilde{\psl2}$.
If a covariant representation $\pi$ is irreducible then there is a 
unique $U_\pi$ satisfying Eq. (\ref{covariance}). 
Hence, in this case, the corresponding generator of rotations
$L_{0}^\pi$ is completely determined by $\pi$. 
Given a covariant representation $\pi$ of $\A$ on a separable Hilbert space 
$\H_\pi$ one has the (isomorphic) inclusions 
$\pi_I(\A(I))\subset \pi_{I^c}(\A(I^c))',\; I\in \I$ \cite{FrG}. Then the
Jones (minimal) index $[ \pi_{I^c}(\A(I^c))' : \pi_I(\A(I)) ]$ is
independent of $I\in \I$ and the {\it statistical dimension} $d(\pi)$ of
$\pi$ is defined by 
\begin{equation}
d(\pi )= [ \pi_{I^c}(\A(I^c))' : \pi_I(\A(I))]^{\frac{1}{2}}.
\end{equation}  
The relation of this definition with the one in \cite{DHR} is given by the 
index-statistics theorem \cite{GuLo96,longo 89-90}.

A {\it conformal subsystem}  of a conformal net
${\cal A}$ is a family 
${\B} =\{ {\B}(I)| \; I\in \I \}$ of nontrivial von Neumann
algebras acting on ${\cal H}_\A$ such that:
\begin{eqnarray}
{\B}(I)\subset {\A}(I) \; \; {\rm for}\; I\in \I;\\
U(\alpha){\B}(I)U(\alpha)^{-1}={\B}(\alpha I)
\; \; {\rm for}\; I , \in \I;\\  
{\B}(I_1)\subset {\B}(I_2)
\; \; {\rm for}\; I_1 \subset I_2, \;\;I_1 , I_2 \in \I.
\end{eqnarray}

We shall use the notation $\B\subset \A$ for conformal subsystems. 
Note that $\B$ is not in general a conformal net since $\Omega$ is not 
cyclic for the algebra $\B(\s1):=\bigvee_{I\in \I} \B (I)$ unless $\B=\A$.
However one gets
a conformal net $\B_0$ by restriction of the algebras $\B(I),\;I\in\I,$
and of the representation $U$ to the closure $\H_\B$ of $\B(\s1)\Omega$.
Since the map $$b\in \B(I) \mapsto b |_{\H_\B} \in \B_0(I) $$
is an isomorphism for every $I\in \I$, we shall, as usual, use the symbol 
$\B$ instead of $\B_0$, specifying, if necessary, when $\B$ acts on
$\H_\A$ or on $\H_\B$.  

Given a conformal subsystem $\B\subset\A$ the index of the subfactor 
$\B(I) \subset \A(I)$ does not depend on $I$ and is denoted $[\A:\B]$. 

\section{Restricting representations} 
We now consider restriction of representations. Given a subsystem
$\B\subset \A$ and a representation $\pi $ of $\A$ one can
define a representation $\pi^{rest}$ by 
\begin{equation}
\pi^{rest}_I=\pi_I|_{\B(I)}\quad I\in\I.
\end{equation}
Then the following holds \cite{Reh2} (cf. also \cite[Section 3]{xu 2000a}). 
We include the proof for the convenience of the reader.  

\begin{proposition}
\label{indrest}
 For every conformal subsystem $\B\subset \A$  and 
covariant representation $\pi$ of $\A$ on a separable Hilbert space 
we have
\begin{equation} 
d(\pi^{rest})=[\A:\B]  d(\pi). 
\end{equation}
\end{proposition}
\begin{proof} For $I\in \I$ we have 
$d(\pi^{rest})^2=[ \pi_{I^c}(\B(I^c))' :\pi_I(\B(I))]$. Consider the inclusions
$$\pi_I(\B(I))\subset\pi_I(\A(I)) \subset \pi_{I^c}(\A(I^c))'\subset 
\pi_{I^c}(\B(I^c))'.$$ 

Then, the multiplicativity of the index \cite{longo 92} implies that
$d(\pi^{rest})^2 $
is equal to
$$
[\pi_{I^c}(\B(I^c))':\pi_{I^c}(\A(I^c))']
[\pi_{I^c}(\A(I^c))':\pi_I(\A(I))][\pi_I(\A(I)):\pi_I(\B(I))]. 
$$
Since $\pi_I$ is an isomorphism for every $I\in\I$ we have 
\begin{align*}
[\pi_{I^c}(\B(I^c))':\pi_{I^c}(\A(I^c))']&=[\pi_{I^c}(\A(I^c)):\pi_{I^c}(\B(I^c))]\\
&=[\A:\B]
\end{align*}
and similarly 
$$[\pi_I(\A(I)):\pi_I(\B(I))]=[\A:\B].$$
It follows that 
$$d(\pi^{rest})^2 = [\A :\B]^2 d(\pi)^2.$$
\end{proof}

\begin{remark}{\rm  If $N\subset M$ is an inclusion of infinite factors
acting on a separable Hilbert space and $\rho$ is a (normal, unital)
endomorphism of $M$ one can define an endomorphism $\rho^{rest}$
of $N$ by 
\begin{equation}
\rho^{rest} := \gamma \circ \rho|_N,
\end{equation}
where $\gamma$ is Longo's canonical endomorphism \cite{longo 89-90}.
As discussed in \cite{LR} the
mapping $\rho \mapsto \rho^{rest}$ (called $\sigma$ restriction in
\cite{BE}) corresponds in a natural way to the restriction of representations of
a net. In fact a similar argument to the one used in the proof of the previous
proposition shows that 
\begin{equation}
d(\rho^{rest})=[M:N]d(\rho).
\end{equation}
 Here the dimension $d(\rho)$ 
of an endomorphism $\rho$ of a factor $M$ is given by the square root of the
index of the subfactor $\rho (M)\subset M$.}
\end{remark}

Let $I_1, I_2\in \I$ have disjoint closures,  let $I_3, I_4\in \I$ be 
the interiors of the connected components of $\s1 /(I_1 \cup I_2)$ and
let $\A$ be a conformal net on $\s1$. 
The inclusion 
\begin{equation} 
\A(I_1)\vee\A(I_2) \subset (\A(I_3)\vee\A(I_4) )'
\end{equation}
is called a {\it 2-interval inclusion}. 
A conformal net $\A$ is said to be {\it completely rational} if it is
split, strongly additive and there is a  2-interval inclusion with finite
index $\mu_\A$ (in this case every  2-interval inclusion has the same
index \cite{KLM}). It has been shown in \cite{KLM}
that a completely
rational net has finitely many sectors which are all covariant with finite
dimension.
Furthermore the following holds
\begin{equation}
\label{muindex}
\mu_\A = \sum_i d(\pi_i)^2,
\end{equation}
where for each sector of $\A$ a representation $\pi_i$ has been chosen. 

We now consider a conformal subsystem $\B$ of a completely rational net 
$\A$ such that the index $[\A :\B ]$ is finite. Then $\B$ is completely 
rational \cite{longo 2001} and the index $\mu_\B$ is given by 
(\cite[Proposition 24.]{KLM})
\begin{equation}
\label{submuindex}
\mu_\B=[\A :\B ]^2 \mu_\A.
\end{equation}
We say that a sector of $\B$ is  {\it untwisted} if it is contained in
$\pi^{rest}$ for some irreducible representation $\pi$ of $\A$ on a
separable Hilbert space. If it is not untwisted we say that it is 
{\it twisted}. 
 
For every sector of $\B$ we choose a corresponding representation 
$\sigma_i$ of $\B$. 
 Let $\U$, ($\T$) be the set of untwisted  (twisted) sectors of $\B$.
We define
\begin{eqnarray}
{\mu^u}_\B = \sum_{[\sigma_i]\in \U} d(\sigma_i)^2, \\
{\mu^t}_\B = \sum_{[\sigma_i]\in \T} d(\sigma_i)^2.
\end{eqnarray}
Clearly ${\mu}_\B ={\mu^u}_\B+{\mu^t}_\B$.
In the case where $\B\subset \A$ is an orbifold inclusion, namely $\B$ is
the fixed points net for the action of a (non trivial) finite group $G$ of
internal symmetries of 
$\A$, it has been shown by Xu \cite{xu 2000b} that the set of twisted
sectors is not empty.
Actually  Proposition \ref{indrest} implies the existence of such sectors
even when there is no underlying group action. 
\begin{proposition}
\label{twisted}
Let $\B$ be a proper conformal subsystem of a completely rational 
net $\A$, with finite index $[\A:\B]$. Then the set of twisted sectors 
of $\B$ is not empty and in fact ${\mu^t}_\B\geq 2$. 
\end{proposition}

\begin{proof} Let $\pi_i,\;i= 0,1,...,n$ be inequivalent
irreducible representations exausting all sectors of $\A$ and let
$\pi_0$ be the vacuum representation. The set $\U$ of untwisted sectors of
$\B$ can be decomposed into disjoint subsets $\U_i,\;i=0, 1,...,n$ in the 
following way:
$\U_0$ is the set of sectors of $\B$ which are contained in 
${\pi_0}^{\it rest}$ and $\U_k,\; k=1,...,n$  is the
the set of sectors contained in ${\pi_k}^{\it rest}$ but not in 
${\pi_i}^{\it rest}, i=0,...,k-1$. It follows from Proposition 
\ref{indrest} and Eq. (\ref{submuindex}) that

\begin{align*}
\sum_i d({\pi_i}^{\it rest})^2 &=[\A:\B]^2 \cdot \sum_i d({\pi_i})^2 \\
&=[\A:\B]^2 \mu_\A \\
& =\mu_\B.
\end{align*}

Therefore

\begin{align*}
{\mu^t}_\B &={\mu}_\B - {\mu^u}_\B \\
&=\sum_i d({\pi_i}^{\it rest})^2 - \sum_{[\sigma_k]\in \U} d(\sigma_k)^2 
\\
&\geq \sum_i (\sum_{ [\sigma_k]\in\U_i}d(\sigma_k))^2 -
\sum_i(\sum_{[\sigma_k]\in \U_i} d(\sigma_k)^2) \\
&\geq  (\sum_{[\sigma_k]\in\U_0}d(\sigma_k))^2 -
\sum_{[\sigma_k]\in \U_0} d(\sigma_k)^2 \geq 2, 
\end{align*}

where the last inequality follows from the fact that $\U_0$ has 
two or more  elements when $\B \neq \A$.
\end{proof}

\section{Virasoro algebra and infinite dimension} 
We begin this section with the following easy consequence of Proposition
\ref{indrest}
\begin{proposition}
\label{instat}
 Let $\B$ be a conformal subsystem of a net $\A$ with
infinite index $[\A:\B]$. Assume that there exists a covariant
representation $\pi$ of $\A$ on a separable Hilbert space whose
restriction to $\B$ is irreducible.
Then $[\pi^{rest}]$ is a covariant sector of $\B$  with infinite
statistical dimension. 
\end{proposition} 

We now come to the sectors of the conformal 
net $\A_{\rm Vir}$ generated by the Virasoro algebra with $c=1$.
 We shall use 
the fact that  $\A_{\rm Vir}$ can be considered has a conformal subsystem 
of the net $\A$ generated by a $\u1$ current $J(z)$. 
The net $\A$ is defined has follows, see \cite{BMT,BS-M} for more
details. The Hilbert space $\H_A$ carries a
strongly continuous unitary representation $U$ of $\psl2$ with positive 
energy and a unique (up to a phase) $U$-invariant unit vector $\Omega$. 
The $\u1$ current $J(z),\; z\in \s1$ is defined as operator valued 
distribution on $\H_\A$. Namely  the operators  
\begin{equation}
J(u) = \int\frac{dz}{2\pi i}J(z)u(z)\;\; u\in C^\infty(\s1)
\end{equation}
have a common invariant dense domain $\D$ containing $\Omega$
which is also $U$-invariant. For each $\psi \in \D$ the mapping
$u\mapsto J(u)\psi$ is linear and continuous from $C^\infty(\s1)$ to
$\H_\A$. Moreover the vacuum $\Omega$ is cyclic for the polynomial algebra 
generated by the smeared currents $J(u),\; u\in C^\infty(\s1)$.
 
The current $J(z)$ satisfies the canonical commutation relations

\begin{equation}
[J(z_1),J(z_2)] = -\delta ' (z_1- z_2), 
\end{equation} 
where the Dirac delta function $\delta (z_1- z_2)$ is defined with respect
to the complex measure $\frac{dz}{2\pi i}$,  the hermiticity condition 
\begin{equation} 
J(z)^*=z^2 J(z),
\end{equation} 
and the covariance property 
\begin{equation} 
U(\alpha)J(u)U(\alpha)^*=J(u_\alpha),\;u\in C^\infty(\s1),
\end{equation}
where $u_\alpha(z):=u(\alpha^{-1}z)$. 
For every real test function $u\in C^\infty(\s1)$ the operator $J(u)$ is 
essentially self-adjoint and the unitaries $W(u):={\rm e}^{iJ(u)} $
satisfy the Weyl relations
\begin{equation}
W(u)W(v)=W(u+v){\rm e}^{-\frac{A(u,v)}{2}}, 
\end{equation}
where $A(u,v):= \int\frac{dz}{2\pi i}u'(z)v(z)$.
For every $I\in \I$ the local von Neumann algebra $\A(I)$ is defined by
\begin{equation}
\A(I)= \{W(u)| u\in C^\infty(\s1)\; {\rm  real},\; {\rm supp}\; u 
\subset I\}''
\end{equation} 
and one can show that the family ${\A(I),\; I\in\I}$ is a conformal net on
$\s1$. 
Next we define the conformal subsystem $\A_{\rm Vir}$ generated by the
Virasoro algebra with central charge $c=1$. 
First consider the (formal) Fourier expansion of the U(1)-current
\begin{equation}
J(z)=\sum_n J_n z^{-n-1},
\end{equation}    
where the Fourier modes $J_n,\;n\in \ZZ$ satisfy 
\begin{eqnarray}
\label{lieu1}
[J_n,J_m]=n\delta_{n+m,0} \\
{J_n}^*=J_{-n}.
\end{eqnarray} 
One can define an energy-momentum tensor $T(z)$ by the Sugawara
construction
\begin{equation}
T(z)=\frac{1}{2}:J(z)^2: = \frac{1}{2}(J_+ (z)J(z)+J(z)J_- (z)),
\end{equation}
where, 
\begin{equation}
J_+(z)=J(z)-J_-(z)=\sum_{n=1}^\infty J_{-n} z^{n-1}.
\end{equation}

The Fourier modes in the expansion 
\begin{equation}
T(z)=\sum_n L_n z^{-n-2}
\end{equation}
satisfy the Virasoro Algebra
\begin{equation}
\label{lievir} 
[L_n,L_m]=(n-m)L_{n+m} +\frac{c}{12}n(n^2-1)\delta_{n+m,0}
\end{equation}
with central charge $c=1$, and the hermiticity condition 
\begin{equation}
{L_n}^*=L_{-n}.
\end{equation}  
According to our previous notations (the closure of) $L_0$ is the
positive self-adjoint generator of the restriction of $U$ to the
one-parameter subgroup of rotations. 

For $f\in C^\infty(\s1)$ the operator 
\begin{equation}
T(f)= \int\frac{dz}{2\pi i}T(z)f(z)
\end{equation}
is well defined on $\D$ and is essentially self-adjoint when 
$z^{-1}f(z)$ is real. The conformal subsystem $\A_{\rm Vir}\subset \A$ is
then defined by
\begin{equation}
\A_{\rm Vir} (I)= \{{\rm e}^{iT(f)}| f\in C^\infty(\s1),\; z^{-1}f(z)\; 
{\rm real},\; {\rm supp}\; f\subset I\}'',\;\;I\in\I.
\end{equation}

Representations of the net $\A$ have been studied in \cite{BMT}. For every 
$q\in \RR$ one can define a covariant irreducible representation
(BMT-automorphism) $\alpha_q$ on $\H_q=\H_\A$ such that 
\begin{equation}
{\alpha_q}_I(W(u))={\rm e}^{ q\int\frac{dz}{2\pi i}z^{-1}u(z)}W(u),
\end{equation}
for $I\in \I$, $u\in C^\infty(\s1)$ with support in $I$. Such
representations have dimension $d(\alpha_q)=1$ and correspond to 
(unitary) positive 
energy representations of the Lie algebra (\ref{lieu1}) with lowest weight 
$q$ \cite{BMT2}. Note that $\alpha_0$ is the vacuum representation of 
$\A$. Analogously to each representation of the Virasoro algebra 
(\ref{lievir}) with
central charge $c=1$ and lowest weight $h\in \RR_{+}$ one can associate
a covariant irreducible representation $\pi_h$ of $\A_{\rm Vir}$ which can
be realized has a subrepresentation of ${\alpha_q}^{rest}$ if 
$h=\frac{1}{2}q^2$, see \cite{BS-M}. The characters of the representations 
$\alpha_q$, $q\in \RR$, are given by (see e.g. \cite[Section 2.2.]{KaRa} )
\begin{equation}
\chi_q (t) ={\rm Tr}( t^{L_0 ^{\alpha_q}})= t^{\frac{1}{2}q^2} p(t)\;
\;t\in (0,1),
\end{equation}
where $p(t)=\prod_{n=1}^{\infty} (1-t^n)^{-1}$. Moreover, for the
representations
$\pi_h$, $h\in \RR_+$ and $t\in (0,1)$, by the results 
in \cite{Kac} the following hold 

\begin{eqnarray}
\chi^h (t):={\rm Tr}( t^{L_0 ^{\pi_h}})=t^{j^2}(1-t^{2|j|+1}) p(t),\;
h=j^2, j\in \frac{1}{2}\ZZ, \\
\chi^h (t):={\rm Tr} (t^{L_0 ^{\pi_h}})=t^h p(t),\; h\neq j^2, j\in
\frac{1}{2}\ZZ.
\end{eqnarray}

\begin{lemma} $[\A : \A_{\rm Vir}]=\infty.$    
\end{lemma}
\begin{proof} As a consequence of  Proposition \ref{indrest} we have 
$[\A : \A_{\rm Vir}]=d({\alpha_0}^{rest})$. Moreover it follows 
from the equality $\chi_0 (t)=\sum_{j=0}^{\infty} \chi^{j^2}(t)$ that 
$${\alpha_0}^{rest}=\oplus_{j=0} ^{\infty} \pi_{j^2}$$
and this implies infinite index.
\end{proof}
\begin{lemma}{\rm (cf. \cite[Theorem 6.2.]{KaRa})} If $h=\frac{1}{2}q^2$,
$q\notin \frac{1}{\sqrt{2}}\ZZ$,
then $\pi^h={\alpha_q}^{rest}$.
\end{lemma}
\begin{proof} If $h=\frac{1}{2}q^2$ $\pi^h$ is a subrepresentation of 
${\alpha_q}^{rest}$ on a $U_{{\alpha_q}}$-invariant subspace
$\H_h\subset \H_q$. Moreover, if 
$q\notin \frac{1}{\sqrt{2}}\ZZ$ then $\chi^h(t)=\chi_q(t)$ and hence
$\H_h = \H_q$. Accordingly we have $\pi^h={\alpha_q}^{rest}$.
\end{proof}
The following theorem is a direct consequence of Proposition 
\ref{instat} and the previous two lemmata. 
\begin{theorem} 
\label{infcont}
If $[\pi_h]$ belongs to the continuum sectors of 
$\A_{\rm Vir}$, i.e. $h\in \RR_+$, $h\neq j^2$, $j\in\frac{1}{2}\ZZ$, 
then it has infinite statistical dimension. 
\end{theorem}
\begin{remark}{\rm It has been shown by Rehren \cite{Reh} that 
if $h=j^2$, $j\in \ZZ$, then $d(\pi_h)=2|j|+1$ and the same formula is
expected to hold for every $j\in \frac{1}{2}\ZZ$.}
\end{remark}

\medskip

\noindent{\bf Acknowledgements.} The author would like to thank Roberto
Longo for some stimulating conversations. He also thanks Roberto Conti and 
Karl-Henning Rehren for useful comments on the manuscript. 

This work was essentially completed when the author was at the 
Dipartimento di Matematica of the Universit\`a di Roma Tre thanks to a 
Post-Doctoral grant of this University.

\end{document}